\documentclass[11pt]{article}
\usepackage{amsmath}
\usepackage{amssymb}
\begin{document}
\title{On configuration of Limit Cycles in certain planner vector  fields}
\author{Ali Taghavi \\ Institute for Advanced Studies in Basic Sciences \\  Zanjan 45195-159, Iran}
\maketitle
Let $X_\lambda$ be a one parameter family of vector field on the
plane satisfying $Det(X_{\lambda_1},X_{\lambda_2})>0$ for
$\lambda_1>\lambda_2$.\\
This means every solution of $X_{\lambda_1}$ is transverse to
solutions of $X_{\lambda_2}$. We call $X_\lambda$ a family  of
rotated vector fields (Note  that such family can be defined on
any symplictic manifold,and each $X_\lambda$ is transvers to
isotropic or lagrangian submanifold invariant under a
$X_{\lambda_0}$ thus it would be interesting to equip a symlictic
manifold to new volume symlictic form, in order to facilitate in
working with a family which is not "rotated family" with respect
to usual symplictic form). This phenomenon have been presented by
Duff $[1]$. In this note, we prove three observation using the
properties of rotated families (In third observation, however, we
do not have a rotated family, but the argument is similar to the
methods in rotated
vector fields).\\
{\bf Proposition i.} The quadratic system
\begin{equation}
\begin{cases}
 \dot{x}=y+ax^2+by^2+cx\\
 \dot{y}=-x+dx^2+fxy
\end{cases}
\end{equation}
can not have two limit cycles with disjoint interiors.\\
{\bf Proposition ii.} The  Lienard system
\begin{equation}
\begin{cases}
 \dot{x}=y+ax^5+bx^3+cx\\
 \dot{y}=-x
\end{cases}
\end{equation}
has  a semistable limit cycle  if  and only if $bc<0$  and
$a=\phi(b,c)$, where $\phi$ is  a unique analytic  function.\\
{\bf Proposition iii.} Lienard system
\begin{equation}
\begin{cases}
 \dot{x}=y+(x^4-2x^2)\\
 \dot{y}=\varepsilon(a-x)
\end{cases}
\end{equation}
has at least one limit cycle if and only if $0<|a|<1$.\\
{\bf Remark.} Proposition 1 could be in particular due to question
posed in [8] about coexistence of two limit cycles with disjoint
interior in quadratic system. Proposition 2 is actually giving a
partial answer to a question about multiplicity of limit cycle in
Rychkov-Lienard system (see[9-page 261], and [6]). Proposition 3
would try suggesting a counterexample of a system
$$
\begin{cases}
 \dot{x}=y-(ax^4+bx^3+cx^2+dx)\\
 \dot{y}=-x
\end{cases}
$$
with at least two limit cycles. See conjecture in [3] about system
(3), it seems that no duck limit cycles could be existed (Due to
intuitions from canard solutions described in [2]). From other
hand Proposition 3 assert that we have at least one limit cycle.
This shows that perhaps for $\varepsilon$ and $a$ small, the limit
cycles bifurcate  from infinity, however the minimum values of
$y$-coordinates of such limit cycle(s) can not be less than -1,
using Remark 3 in [7]. Thus it would be interesting investigation
of the number of limit cycles of (3) or adding a term $\varepsilon
x^3$ to first line of (3). I thank Professor Roussarie that he
explained about canard solutions and suggested the latest system
as a possibly candidate for counterexample to Pugh's conjecture.\\\\
{\bf Proof  of Proposition 1.} This is proved in three steps\\
{\bf i.} If a limit  cycle  surround the origin then
$cd(2a+f)>0$,\\
{\bf ii.} If  a limit  cycle  does  not intersect the  line x=0
and has positive (negative) orintion    then $cd(2a+f)<0,(>0)$,\\
{\bf iii.} If $cd(2a+f)=0$  then two  limit cycles  can  not
coexist.\\
Assume that all 3 statements in above  are  proved,let $\gamma_1$
,$\gamma_2$ are  two  limit  cycles  with  disjoint interiors, by
{\bf ii} and {\bf iii} at least one of the $\gamma_1$ and
$\gamma_2$ must intersect the y-axis and we may assume that the
origin lies in $\gamma_1$ (for if not we translate the singularity
inside of $\gamma_1$  to origin. From {\bf i} and {\bf ii} we
obtain that $\gamma_2$ must also intersect the line $x=0$.
Therefore Both $\gamma$'s do not intersect the line $-1+dx+fy=0$
because any closed orbit of a quadratic system can surround only
one singularity[10]. Now we add $-cx$ to  first equation of (1)
and we obtain a limit cycle for (1) when $c=0$, while is
impossible, see [9-page 315].\\
{\bf Proof  of step iii.} when $c=0$, (1) does not have a limit
cycle because of the reason mentioned in above two line, if
$2a+f=0$, divergence of (1) is constant thus there is non limit
cycle and if $d=0$, we have at most one limit cycle, see [9], in
which the origin does not lie because for $c=0$ and $d=0$ the
origin is a center: (note  that in a rotated vector field family,
if we have a center for a parameter $\lambda_0$ we could not have
limit cycle for other values of $\lambda$.\\
{\bf Proof  of step i.} If the origin lies inside a limit cycle
then $cd(2a+f)$ is not $0$ and if it is negative we add $-cx$ to
first equation and a contradiction is follows.\\
{\bf Proof of step ii.} Note that if a limit cycle does not
intersect the $y$-axis,then $x$ values of its point has the same
sign as the $\frac{-c}{2a+f}$ and by computation of
$\int_{\gamma}(-1+dx+fy)dy$ we find that it has the same sign as
$d$ (for positive orient of parameterizations of limit cycle
$\gamma$, then $cd(2a+f)$ is negative. Similar consideration hold
for negative orient and the proof is completed.\\
{\bf Proof  Of  Proposition ii.} It is proved in [6] that system
(2) has at most two limit cycles. In fact this result is true
counting multiplicity: Let $P(y)$  be the poincare  map defined on
positive $y$-axis. Then $P'(y)=\frac{y}{P(y)}e^{h(y)}$ where
$h(y)=\int_0^{T(y)}\text{divergence of(2)}$, $T(y)$ is the time of
first return.\\
Assume that $P(y)=y_0$ and $P'(y_0)=1$, the computation in [6],
actually shows that $h'(y_0)\neq 0$ so $P''(y_0)\neq 0$ then (2)
has at most two limit cycle counting multiplicity. Now We present
a global bifurcation diagram of semi-stable limit cycle for (2).
If $bc>0$ then by lienard theorem [5], there is no semistable
limit cycle. Assume that $bc<0$. For $a=0$, system (2) has  a
unique hyperbolic limit cycle. We can assume $c<0$ and $b>0$, if
$a<0$ and $|a|\lll1$, then another limit cycle would born at
infinity. If for some $a_0<0$, two limit cycles exist, then the
same holds for $a_0<a<0$, because if $\gamma_1$ and $\gamma_2$
would be two limit cycles for $(2)_{a_0}$, then both of $\gamma_1$
and $\gamma_2$ are closed curve without contact for $(2)_a$ for
all $a_0<a<0$. Now Compare the direction of $(2)_a$ on $\gamma_1$
and $\gamma_2$ with stability of origin and infinity. On the other
hand for fixed $c<0$ , $b>0$, if $|a|$ is sufficiently large
$(a<0)$, then the derivative of energy  does not change sign.
Therefore there exist a unique $a_0=\phi(b,c)$ such that (2) has a
semistable limit cycle. $a_0$ is unique because from any
semiustable limit cycle, two limit cycles could be created. Now
all conditions of Theorem 2 in [4] satisfy and proposition 2 is
proved.\\
{\bf Proof of proposition iii.} For $|a|\geq 1$ there is no limit
cycle using proposition in [3], after change of coordinate
$x:=x+a , y:=y+a^4-2a^2$.\\
For $a=0$ the system (3) has a center whose region of closed
orbits is bounded by a unique orbit $\gamma$ asymptotic to the
graph of $y=x^4-2x^2$ and  $\gamma$ is below this graph, thus
$\gamma$ is a curve without contact for $(3)_a$ and the
singularity is attractive. Thus Poincare Bendixon theorem convert
to existence of at least one limit cycles.\\
{\bf  Remark.} It Is obvious that a multiple limit cycle (with
arbitrary finite large multiplicity) can produce at most two limit
cycles with one parameter perturbation in a rotated family. How
much this results remain valid in the case of infinite
multiplicity? See [5-page 387].

\end{document}